\documentclass[12pt]{amsart}
\usepackage{amssymb}
\usepackage{amsmath}
\newcommand\al{\alpha}

\newcommand\Aa{\mathcal{A}}
\newcommand\AND{\quad\text{and}\quad}
\newcommand\Bb{\mathsf{B}}
\newcommand\C{\mathbb C}
\newcommand\Cc{\mathsf{C}}
\newcommand\de{\delta}
\newcommand\ep{\varepsilon}
\newcommand\Ex{\mathsf{E}}
\newcommand\Gp{\mathfrak{G}}
\newcommand\lb{\boldsymbol{\ell}}
\newcommand\Ll{\mathsf{L}}
\newcommand\Lp{\mathfrak{L}}
\newcommand\Mm{\mathcal{M}}
\newcommand\N{\mathbb N}
\newcommand\Prob{\mathsf{Pr}}
\newcommand\R{\mathbb R}
\newcommand\rb{\mathbf{r}}
\newcommand\Ss{\mathsf{S}}
\newcommand\Sf{\mathfrak{S}}
\newcommand\supp{\operatorname{\sf supp}}
\newcommand\tb{\mathbf{t}}
\newcommand\Uu{\mathcal{U}}
\newcommand\uno{\mathbf{1}}
\newcommand\wh{\widehat}
\newcommand\wt{\widetilde}
\newcommand\Xx{\mathsf{X}}
\newcommand\Z{\mathbb Z}

\numberwithin{equation}{section}

\newtheoremstyle{mythm}
  {9pt}
  {9pt}
  {\itshape}
  {0pt}
  {\bfseries}
  {}
  { }
  {\thmnumber{(#2)}\thmname{ #1}\thmnote{ #3}}

\newtheoremstyle{mydef}
  {9pt}
  {9pt}
  {\normalfont}
  {0pt}
  {\bfseries}
  {}
  { }
  {\thmnumber{(#2)}\thmname{ #1}\thmnote{ #3}}

\theoremstyle{mythm}
\newtheorem{thm}[equation]{Theorem.}
\newtheorem{pro}[equation]{Proposition.}
\newtheorem{lem}[equation]{Lemma.}
\newtheorem{cor}[equation]{Corollary.}

\theoremstyle{mydef}
\newtheorem{dfn}[equation]{Definition.}
\newtheorem{exa}[equation]{Example.}

\begin{document}$\,$ \vspace{-1truecm}
\title{On recurrence of reflected \\ random walk on the half-line}
\author{\bf Marc PEIGN\'E and Wolfgang WOESS \\ \ \\
With an appendix on results  of Martin BENDA}
\address{\parbox{.8\linewidth}{Laboratoire de Math\'ematiques
et Physique Th\'eorique\\
Universit\'e Francois-Rabelais Tours\\
F\'ed\'ration Denis Poisson -- CNRS\\
Parc de Grandmont, 37200 Tours, France\\}}
\email{marc.peigne@univ-tours.fr}
\address{\parbox{.8\linewidth}{Institut f\"ur Mathematische Strukturtheorie 
(Math C),\\ 
Technische Universit\"at Graz,\\
Steyrergasse 30, A-8010 Graz, Austria\\}}
\email{woess@TUGraz.at}
\date{\today} 
\thanks{The second author acknowledges support by a visiting professorship
at Universit\'e de Tours}
\subjclass[2000] {60G50; 
                  60J05 
                  }
\keywords{reflected random walk, recurrence, invariant measure, process of
reflections, local contractivity, stochastic iterated function system}
\begin{abstract}
Let $(Y_n)$ be a sequence of i.i.d. real valued random variables.
Reflected random walk $(X_n)$ is defined recursively by $X_0=x \ge 0$,
$X_{n+1} = |X_n - Y_{n+1}|$. In this note, we study recurrence of this
process, extending a previous criterion. This is obtained by determining
an invariant measure of the embedded process of reflections.
\end{abstract}

\maketitle

\markboth{{\sf M. Peign\'e and W. Woess}}
{{\sf Reflected random walk}}
\baselineskip 15pt

\section{Introduction}\label{intro}

Reflected random walk was described and studied by {\sc Feller~\cite{Fe}}; 
apparently, it was first considered by {\sc von Schelling~\cite{Sch}} 
in the context of telephone networks.
 
Let $(Y_n)_{n \ge 0}$ be a sequence of i.i.d. real valued random variables,
and let $S_n = Y_1 + \ldots + Y_n$ be the classical associated random
walk. Reflected random walk is obtained by considering a non-negative
initial random variable $X_0$ independent of the $Y_n$ and considering 
$X_0 - S_n$, $n =0, 1, ...$,  as long as this is non-negative. When it
becomes negative, we change sign and continue from the new (reflected)
point by subtracting $Y_{n+1}$, $Y_{n+2}$, ..., until the next reflection, 
and so on. Thus, we consider the Markov chain $X_n$ given by
$X_{n+1} = |X_n - Y_{n+1}|$. We are interested in recurrence of this process
on its essential (i.e., maximal irreducible) classes. 

We start by considering the situation when $Y_n \ge 0$ (of course excluding 
the trivial case $Y_n \equiv 0$), so that the increments of $(X_n)$ are
non-positive except possibly at the moments of reflection. In this case,
{\sc Feller~\cite{Fe}} and {\sc Knight~\cite{Kn}} have computed an invariant measure 
for the process when the $Y_n$ are non-lattice random variables, while 
{\sc Boudiba~\cite{Bo1}}, \cite{Bo2} has provided such a measure
when the $Y_n$ are lattice variables. {\sc Leguesdron~\cite{Le}}, 
{\sc Boudiba~\cite{Bo2}} and {\sc Benda~\cite{Be}} have also studied its 
uniqueness (up to constant factors). When that invariant measure has
finite total mass -- which holds if and only if $\Ex(Y_1) < \infty$ -- the 
process is (topologically) recurrent: with probability $1$, 
it returns infinitely often to each open set that is charged by the invariant
measure.  

Our main result is that reflected random walk is still 
recurrent when $Y_n \ge 0$ and 
$\int_0^{\infty} \Prob[Y_1 \ge t]^2 \,dt < \infty\,$; see \S 3 for the
case when the $Y_n$ are lattice random variables, and \S 4 for the non-lattice
case.
The result is based on considering the process of reflections, that is,
reflected random walk observed at the instances of reflection, see \S 2. 
We determine an invariant measure for the latter. The above ``quadratic tail'' 
condition holds if and only if that measure is finite. This holds, in 
particular, when $\Ex(Y_1^{1/2}) < \infty$. 

Subsequently, in \S 5, we also consider the case when the $Y_n$ may assume 
negative as well as positive values. Reflected random walk is of interest when 
$\limsup_n S_n = \infty$ almost surely. Let $Y_1 = Y_1^+ - Y_1^-$ be the 
decomposition into positive and negative part. If 
$\Ex(Y_1^-) < \Ex(Y_1^+)$ then the situation is similar to the case when
$Y_1 \ge 0$ a.s., and we get recurrence when  
$\Ex\bigl(\sqrt{Y_1^+}\,\bigr) < \infty\,$. 
If the $Y_n$ are centered, that is, $0 < \Ex(Y_1^-) = \Ex(Y_1^+)$,
then we get recurrence under the moment condition
$\Ex\Bigl(\sqrt{Y_1^+}^{\,3}\Bigr) < \infty\,$, which turns out to be almost
sharp. 

Our methods are based on interesting and useful work of {\sc M. Benda}
in his PhD thesis \cite{Be} (in German)  and the two subsequent preprints 
\cite{Be1}, \cite{Be2} which have remained unpublished. For this reason,
we outline those results in the Appendix (\S 6). 

\section{The process of reflections}\label{sect:reflection}

In this and the next two sections, we suppose always that $(Y_n)$ is a sequence of 
i.i.d, non-constant, non-negative random variables. Let $\mu$ be the
(common) distribution of $Y_n$, a non-degenerate probability 
measure on $[0\,,\,\infty)\,$, and $F(x) = F_{\mu}(x) = \mu([0\,,\,x])$ the
associated distribution function ($x \ge 0$). Denote by $\mu^{(n)}$ its 
$n$-th convolution
power, the distribution of $S_n$, with $\mu^{(0)} = \de_0$. 
Since $S_n \to \infty$ almost surely, the potential
\begin{equation}\label{eq:potential}
\Uu = \sum_{n=0}^{\infty} \mu^{(n)}
\end{equation} 
defines a Radon measure on $[0\,,\,\infty)\,$, that is, $\Uu(B) < \infty$
if $B$ is a bounded Borel set. 

Now consider the sequence of stopping times $(\rb(k))_{k \ge 0}$, where 
$\rb(0) = 0$, and $\rb(k)$ ($k > 0$) is the time of the $k$-th reflection:
\begin{equation}\label{eq:reflectiontime}
\begin{aligned}
\rb(k+1) &= \inf \bigl\{ n > \rb(k) : X_n = -(X_{n-1} - Y_n) \bigr\} \\
&=  \inf \bigl\{ n > \rb(k) :  (Y_{\rb(k)+1} + \ldots + Y_{n-1} + Y_n)
\ge X_{\rb(k)} \bigr\}\,.
\end{aligned}
\end{equation} 
Once more because $S_n \to \infty$, each $\rb(k)$ is finite almost surely. 
We call the embedded process $R_k = X_{\rb(k)}$ the 
\emph{process of reflections.}

\begin{lem}\label{lem:reflection} The process of reflections is a Markov
chain with transition probabilities given as follows: if $B \subset  [0\,,\,\infty)$
is a Borel set, then
$$
q(0,B) = \mu(B) \AND q(x,B)= \int_{[0\,,\,x)} \mu(B+x-w)\, \Uu(dw)\,,\quad
\text{if}\;\; x > 0\,.
$$
\end{lem}

\begin{proof} 
It is clear that $(R_k)$ is a (time-homogeneous) Markov chain. We compute
$$
\begin{aligned}
q(x,B) &= \Prob[R_1 \in B \mid R_0 = x]
= \sum_{n=1}^{\infty} \Prob[ \rb(1) = n\,,\; S_n - x \in B ]\\ 
&= \sum_{n=1}^{\infty} \Prob[ S_{n-1} < x \,,\; S_n - x \in B ]
= \sum_{n=1}^{\infty} \int_{[0\,,\,x)} \Prob[Y_n+w-x \in B]\,\mu^{(n-1)}(dw)\\
&= \int_{[0\,,\,x)} \mu(B+x-w)\, \Uu(dw)\,,
\end{aligned}
$$
as proposed. 
\end{proof}

It is an instructive exercise, relying on the fact that 
$\supp(\mu) \subset [0\,,\,\infty)$, 
to show directly that $q(\cdot,\cdot)$ is stochastic.

Now the idea is the following: if the embedded process of reflections is 
recurrent, then also the original reflected Markov chain must be recurrent.

\section{The lattice case}\label{sect:lattice}\medskip

We start with the discrete case, which is instructive and has to 
be treated separately anyway. Here we suppose that there is $\kappa > 0$
such that $\supp(\mu) \subset \kappa \cdot \N_0$, and we may assume
without loss of generality that $\kappa =1$. (By $\N_0$ we denote the
non-negative integers.) 

The one-step transition probabilities of $(X_n)$
are
\begin{equation}\label{eq:transprob} 
p(x,y) = \begin{cases} \mu(x)\,,&\text{if}\; y=0\,,\\
                       \mu(x+y)\,,&\text{if}\; x < y \,,\\
		       \mu(x-y)+\mu(x+y)\,,&\text{if}\; x \ge y > 0\,.
	 \end{cases}
\end{equation}
We write $p^{(n)}(x,y) = \Prob[X_n=y \mid X_0=x]$ for the $n$-step transition 
probabilities. Set 
$$
d = \gcd \supp(\mu) \AND N = \sup \supp(\mu)\,.
$$
If the reflected Markov chain starts in a deterministic point 
$X_0=x_0 \in [0\,,\,\infty)$, then $(X_n)$ evolves within the state space
$$
\Ss(x_0)=\{ kd \pm x_0 : k \in \Z\} \cap [0\,,\,\infty)\,.
$$ 
Recall that an essential class of a denumberable Markov chain is a subset 
$\Cc$ of the state space which is irreducible and absorbing: 
if $x \in \Cc$ then $p^{(n)}(x,y) > 0$
for some $n$ if and only if $y \in \Cc$. 
The next lemma follows from \cite{Bo2} when the starting point $x_0$ is 
rational, and when it is irrational, it is immediately seen to be true as well.

\begin{lem}\label{lem:essential} The reflected random walk $(X_n)$ starting at
$x_0$  is absorbed after finitely many steps by the essential class
$\Cc(x_0) = \Ss(x_0) \cap [0\,,\,N]\,.$
\end{lem} 

When we speak of recurrence of $(X_n)$ with starting point $x_0$
then we mean recurrence on $\Cc(x_0)$. This is known to be independent of
$x_0$ \cite{Bo2}.

If $N=\infty$ then $\Cc(x_0) = \Ss(x_0)$.
Also, if $\supp(\mu)$ is finite then $\Cc(x_0)$ is finite and
carries a unique invariant probability measure.   
An invariant measure $\nu$ (not necessarily with finite total mass) 
exists always. Its formula is due to
\cite{Bo1}, where only $x_0 \in \Z$ is considered, but it can be adapted to 
the present situation with arbitrary starting point as follows. Set
\begin{equation}\label{eq:invmeas}
\nu(0) = \frac{1-\mu(0)}{2} \AND
\nu(x) = \frac{\mu(x)}{2} + \mu\bigl((x\,,\,\infty)\bigr)\,,\quad\text{if}
\;\;x>0\,.
\end{equation}
Here, we mean of course $\mu(x)= \mu(\{x\})$, so that
$\mu(x)=0$ when $x \in [0\,,\,\infty) \setminus \N_0\,$.
Then the invariant measure $\nu_{x_0}$ on $\Cc(x_0)$ is given by the
restriction of $\nu$ to that essential class: if $B \subset \Cc(x_0)$ then
$\nu_{x_0}(B) = \sum_{x \in B} \nu(x)\,.$ 

\begin{cor}\label{cor:posrec}
The reflected random walk starting at $x_0$ is positive recurrent on
$\Cc(x_0)$ if and only if the first moment $\sum_n n\,\mu(n)$ of $Y_k$ 
is finite. 
\end{cor}

If the reflected random walk is (positive or null) recurrent on $\Cc(x_0)$,
then it follows of course from the basic theory of denumerable Markov chains
that $\nu_{x_0}$ is the unique invariant measure (up to multiplication with
constants). 

We now consider the process of reflections.

\begin{lem}\label{lem:essrefl}
The set $\Cc(x_0)$ is also the unique essential class for $(R_k)$ starting
at $x_0$.
\end{lem}

\begin{proof} Since $\Cc(x_0)$ is the only essential class for $(X_n)$,
we only need to verify that it is an irreducible class for $(R_k)$. We have to 
show that for $x, y \in \Cc(x_0)$, it occurs with positive probability
that $(X_n)$, starting at $x$, reaches $y$ at some reflection time
$\rb(k)$.

There is $m \in \supp{\mu}$ such that $m \ge y$. Then also 
$m-y \in \Cc(x_0)$, and there is $n$ such that $p^{(n)}(x,m-y) > 0$.
But from $m-y$, the reflected random walk can reach $y$ (the reflection
of $-y$) in a single step with positive probability $\mu(m)$, and this
occurs at a reflection time.
\end{proof} 

Our simple new contribution is the following.

\begin{thm}\label{thm:ref-inv}
Set $\;\;\rho(0) = \dfrac{1-\mu(0)}{2}\;\;$ and
$$
\rho(x) = \sum_{k=1}^{\infty} \left(\frac{\mu(x)}{2} + 
\mu\bigl((x\,,\,x+k)\bigr) + \frac{\mu(x+k)}{2}\right) \mu(k)\,,\quad\text{if}
\;\;x>0\,.
$$
Then the restriction $\rho_{x_0}$ of $\rho$ to $\Cc(x_0)$ is 
an invariant measure for the process of reflections $(R_k)$ on $\Cc(x_0)$.
It is unique (up to multiplication by a constant), if
$\nu_{x_0}$ is the unique invariant measure (up to multiplication by a constant) 
for the reflected random walk $(X_n)$ on $\Cc(x_0)$.
\end{thm}

\begin{proof} We first show that $\rho_{x_0}$ is invariant. 
The index $x_0$ will be ommitted whenever this does not obscure
the arguments. Also, note that by its definition, $\rho \equiv 0$ on 
$\Ss(x_0) \setminus \Cc(x_0)$, so that we can think of $\rho_{x_0}$ as a
measure on the whole of $\Ss(x_0)$ with no mass outside $\Cc(x_0)\,$. 

Consider the signed measure  $\Aa$ defined by 
$\Aa(x) = \de_0(x) - \mu(x)$ for $x \ge 0$.
Then we have the convolution formula $\Aa*\Uu = \Uu*\Aa = \de_0$, that is
\begin{equation}\label{eq:convo1}
\sum_{j=0}^n \Aa(j)\,\Uu(n-j) = \de_0(n)\,.
\end{equation}
Now we verify that for each real $x \in (0\,,\,N]\,$, 
\begin{equation}\label{eq:convo2}
\rho(x) = \sum_{k=0}^{\infty} \Aa(k)\,\nu(x+k)\,.
\end{equation}
Indeed, 
the last sum is equal to
$$
\bigl(1 - \mu(0)\bigr) \nu(x) - \sum_{k=1}^{\infty} \mu(k)\,\nu(x+k)
= \sum_{k=1}^{\infty} \mu(k) \,\bigl(\nu(x) -\nu(x+k)\bigr)\,,
$$
which is equal to $\rho(x)$. We remark here that the sum in \eqref{eq:convo2}
is absolutely convergent, since $\nu(\cdot) \le 1$.
Combining \eqref{eq:convo2} with the
inversion formula \eqref{eq:convo1}, we get
$$
\begin{aligned}
\sum_{k=0}^{\infty} \Uu(k)\,\rho(x+k)
&= \sum_{k=0}^{\infty} \Uu(k)\,\sum_{l=0}^{\infty}\Aa(l)\,\nu(x+k+l)\\ 
&= \sum_{n=0}^{\infty}\nu(x+n)\sum_{k=0}^{n} \Uu(k)\,\Aa(n-k)\\ 
&= \sum_{n=0}^{\infty}\nu(x+n)\,\delta_0(n)\,,
\end{aligned}
$$
that is,
\begin{equation}\label{eq:convo3}
\nu(x) = \sum_{k=0}^{\infty} \Uu(k)\,\rho(x+k)\,, \quad x > 0\,.
\end{equation}
If $\sigma$ is any measure on $\Cc(x_0)$ then we write 
$$
\Ex_{\sigma}(\cdot) = \sum_{w \in \Cc(x_0)} \sigma(w) \Ex_w(\cdot)\,,
$$
where $\Ex_w(\cdot)$ denotes expectation when the starting
point is $X_0=w$. We claim that
\begin{equation}\label{eq:expandnu} 
\nu(x) = \Ex_{\rho_{x_0}} \left( \sum_{j=0}^{\rb(1)-1} \uno_x(X_j) \right)\,,
\quad \text{if}\;\; x \in \Cc(x_0)\,.
\end{equation}
Indeed, if $x=0$ then the right hand side of \eqref{eq:expandnu} is
$\rho(0) = \nu(0)$, since the reflected random walk can reach the state $0$
before the first reflection only when it starts at $0$, in which case 
$\rb(1)=1$. If $x > 0$, $x \in \Cc(x_0)$ then the reflected walk starting from 
$w \in \Cc(x_0)$ can reach $x$ before the first reflection only if
$w = x+k$ for some $k \in \N_0$ such that $k = S_j$ 
for some $j \ge 0$. We compute
$$
\begin{aligned}
\Ex_{x+k}\left( \sum_{j=0}^{\rb(1)-1} \uno_x(X_j) \right) 
&= \Ex_{x+k}\left(\sum_{n=1}^{\infty} \uno_n\bigl(\rb(1)\bigr) 
                  \sum_{j=0}^{n-1} \uno_x(X_j) \right)\\ 
&= \sum_{j=0}^{\infty} \Prob[X_j=x, \rb(1) > j \mid X_0 = x+k]\\
&= \sum_{j=0}^{\infty} \Prob[x = x+k - S_j] = \Uu(k)\,.
\end{aligned}
$$
Therefore
$$
\Ex_{\rho_{x_0}} \left( \sum_{j=0}^{\rb(1)-1} \uno_x(X_j) \right)
= \sum_{k=0}^{\infty} \rho(x+k) \,\Uu(k) = \nu(x)\,,\quad\text{if}\;\;x > 0\,,
$$
as proposed. From \eqref{eq:expandnu}, we infer that
$$
\sum_w \nu(w)\,p(w,x) = 
\Ex_{\rho_{x_0}} \left( \sum_{j=1}^{\rb(1)} \uno_x(X_j) \right)\,.
$$
Now $\nu$ satisfies $\sum_w \nu(w)\,p(w,x) = \nu(x)$, and applying 
\eqref{eq:expandnu} once more, we obtain
$$
\Ex_{\rho_{x_0}} \Bigl( \uno_x(X_0) \Bigr)
= \Ex_{\rho_{x_0}} \Bigl( \uno_x(X_{\rb}) \Bigr)\,.
$$
The left hand side is $\rho(x)$, while the right hand side is
$\sum_w \rho(w)q(w,x)$, where $q(\cdot,\cdot)$ is the transition
kernel of the process of reflections. Thus, $\rho_{x_0}$ is invariant
for $(R_k)$ on the state space $\Cc(x_0)$.

\smallskip

We now prove uniqueness. In view of Lemma \ref{lem:essrefl}, this is of course 
obvious by the basic theory of denumerable Markov chains, when 
$\rho_{x_0}(\Cc(x_0)) < \infty$, but this is not supposed in our statement.
 
So let $\bar\rho $ be another invariant measure for $(R_k)$ on $\Cc(x_0)\,$,
again considered on $\Ss(x_0)$ with zero mass outside $\Cc(x_0)\,$.
Using the formula of Lemma \ref{lem:reflection} for the transition probabilities
of $(R_k)$, we get for $y \in \Cc(x_0)$ 
$$
\bar\rho (y) = \sum_{w \in \Cc(x_0)} \bar\rho (w) 
             \sum_{k \in \N_0 : 0 \le k < w} \Uu(k)\,\mu(w+y-k)
= \sum_{k=0}^{\infty}\sum_{w \in \Cc(x_0) : w > k}  \Uu(k)\, \bar\rho (w)\,\mu(w+y-k)
$$
To have a non-zero contribution in the last double sum, $w+y$ has to be 
integer,  $d$ must divide both $k$ and $w+y$, and $x=w-k \in \Cc(x_0)\,$. 
Therefore we can rewrite
$$
\bar\rho (y) 
= \sum_{k=0}^{\infty} \sum_{x \in \Cc(x_0) : x > 0}  
\Uu(k)\,\bar\rho (x+k)\,\mu(x+y)\,.
$$
Now let $x \in \Cc(x_0)\,$, $x > 0$. Again, there is $m \in \supp(\mu)$
with $x \le m$, and $y=m-x \in \Cc(x_0)\,$. Therefore
$$
\sum_{k=0}^{\infty} \Uu(k)\,\bar\rho (x+k) \le \frac{\bar\rho (y)}{\mu(m)} < \infty
$$
for each $x \in \Cc(x_0)$ with $x > 0$. This allows us to define a new measure
$\bar\nu $ on $\Cc(x_0)$ by
$\bar\nu (0) = \bar\rho (0)$, if $0 \in \Cc(x_0)$, and
$$
\bar\nu (x) =  \sum_{k=0}^{\infty} \Uu(k)\,\bar\rho (x+k)\,,\quad\text{if}\;\;x > 0\,,
$$
and a straightforward exercise shows that it is legitimate to apply
the inversion formula \eqref{eq:convo1} to deduce that
$$
\bar\rho (x) =  \sum_{k=0}^{\infty} \Aa(k)\,\bar\nu (x+k)\,,\quad\text{if}\;\;x > 0\,,
$$
The same computations as that lead to \eqref{eq:convo3} and \eqref{eq:expandnu}
show that 
$$
\bar\nu (x) = \Ex_{\bar\rho _{x_0}} \left( \sum_{j=0}^{\rb(1)-1} \uno_x(X_j) \right)
$$
is an invariant measure for $(X_n)$ on $\Cc(x_0)$. By uniqueness of the latter,
$\bar\nu  = c\cdot \nu_{x_0}$ for some $c > 0$. Therefore 
$\bar\rho  = c \cdot \rho_{x_0}$.
\end{proof}

\begin{cor}\label{cor:quadratic}
The total mass of $\rho_{x_0}$ is finite for some (equivalently, every)
starting point $x_0$ if and only if
\begin{equation}\label{eq:quadratic}
\sum_{k=0}^{\infty} \bigl( 1 - F_{\mu}(k) \bigr)^2 < \infty\,.
\end{equation}
\end{cor}

\begin{proof}
We write $H(x) = 1 - F_{\mu}(x)$.
For real $\al \ge 0$, let $\Sigma(\al) = \sum_{k=0}^{\infty} \rho(\al+kd)$.
Let $\al_0$ be the unique number in $(0\,,\,d]$ such that $x_0-\al_0$ 
is an integer multiple of $d$.  If $\al_0=d$ or $\al_0=d/2$ we have 
$\rho\bigl(\Cc(x_0)\bigr) = \rho(0) \de_0\bigl(\Cc(x_0)\bigr) +
\Sigma(\al_0)\,$,
while otherwise $\rho(\Cc(x_0)) = \Sigma(\al_0)+\Sigma(d-\al_0)$. 
%
%
Thus, we prove that for any $\al\in (0\,,\,d]$, we have $\Sigma(\al) < \infty$
if and only if \eqref{eq:quadratic} holds. Recalling that
$\mu(x)=0$ if $x$ is not a multiple of $d$, we compute 
$\Sigma(\al) = \Sigma_0(\al) + \Sigma_1(\al)$, where
$$
\Sigma_0(\al)= \sum_{k=0}^{\infty} \sum_{m=1}^{\infty} 
          \frac{\mu(\al+kd) - \mu(\al+kd+md)}{2}\, \mu(md)
$$
is always finite, and
$$
\begin{aligned} 
\Sigma_1(\al) &= \sum_{k=0}^{\infty} \sum_{m=1}^{\infty} 
          \Bigl( H(\al+kd) - H(\al+kd+md)\Bigr) \, \mu(md)\\
&=\sum_{m=1}^{\infty} \sum_{k=0}^{m-1} H(\al + kd) \, \mu(md)\\
&=\sum_{k=0}^{\infty} H(\al+kd) \sum_{m=k+1}^{\infty} \mu(md)
= \sum_{k=0}^{\infty} H(\al+kd) \,H(kd)\,.
\end{aligned}
$$
Since $H(\cdot)$ is decreasing, on one hand
$$
\Sigma_1(\al) \le \sum_{k=0}^{\infty} H(kd)^2 = 
\frac{1}{d} \sum_{k=0}^{\infty} \bigl( 1 - F_{\mu}(k) \bigr)^2\,,
$$ 	  
and on the other hand 
$$
\Sigma_1(\al) \ge \sum_{k=0}^{\infty} H\bigl((k+1)d\bigr)^2	  
= \frac{1}{d} \sum_{k=0}^{\infty} \bigl( 1 - F_{\mu}(k) \bigr)^2 - H(0)^2\,.
$$
Thus, $\Sigma_1(\al)$ and the sum in \eqref{eq:quadratic} are finite,
resp. infinite, simultaneously.
\end{proof}

The following is now immediate. 

\begin{thm}\label{thm:recurrent}
Suppose that the ``quadratic tail'' condition \eqref{eq:quadratic}
holds. Then the process of reflections $(R_k)$ is positive recurrent
on $\Cc(x_0)$ for each starting point $x_0 \ge 0$. 
If in addition $\Ex(Y_1) = \sum_{k \ge 0} k\,\mu(k) < \infty$,
then the reflected random walk $(X_n)$ is also positive recurrent
on $\Cc_{x_0}$, while it is null recurrent when $\Ex(Y_1)=\infty$.
\end{thm}

Finally, it is easy to relate the ``quadratic tail'' condition 
with a moment condition.

\begin{lem}\label{lem:momenthalf} 
If $\Ex\bigl(\sqrt{Y_1}\bigr) = \sum_{k\ge 0} \sqrt{k}\,\mu(k) < \infty\,$,
then \eqref{eq:quadratic} holds.
\end{lem}

\begin{proof} We use the Cauchy-Schwarz inequality:
$$
\left(\sum_{k=n+1}^{\infty} \mu(k)\right)^2
\le \left(\sum_{k=n+1}^{\infty} \mu(k)\,\sqrt{k}\right)
\left(\sum_{k=n+1}^{\infty} \mu(k)/\sqrt{k}\right)
\le \Ex\bigl(\sqrt{Y_1}\,\bigr)\sum_{k=n+1}^{\infty} \mu(k)/\sqrt{k}\,.
$$
Therefore, 
$$
\sum_{n=0}^{\infty}\bigl( 1 - F_{\mu}(n) \bigr)^2 
\le \Ex\bigl(\sqrt{Y_1}\,\bigr)
\sum_{n=0}^{\infty}\sum_{k=n+1}^{\infty} \mu(k)/\sqrt{k}
= \Bigl(\Ex\bigl(\sqrt{Y_1}\,\bigr)\Bigr)^2\,,
$$
which is finite. 
\end{proof}

\section{The non-lattice case}\label{sect:non-lattice}\medskip

We now consider the case when $\supp(\mu) \subset [0\,,\,\infty)$, but
there is no $\kappa > 0$ such that
$\supp(\mu) \subset \kappa\cdot\N_0$. Again, denote $N = \sup \supp(\mu)$,
and set $\Cc = [0\,,\,N]$ if $N < \infty$, resp. $\Cc=[0\,,\,\infty)$,
if $N = \infty$. The transition probabilities of the reflected
random walk are
$$
p(x,B) = \mu(\{ y : |x-y| \in B \})\,,
$$
where $B \subset [0\,,\,\infty)$ is a Borel set. For the following, we need to
specify in more detail the probability space on which we are working.
This is the product space 
$(\Omega,\Prob) = \Cc^{\N},\mu^{\N}\bigr)$, where $Y_n$ is the
$n$-th projection. It will be convenient to write $X_n^x$ for the reflected 
walk starting at $x \ge 0$, so that $X_0^x=x$ and 
$X_{n+1}^x = |X_n^x - Y_{n+1}|$ as in the 
Introduction. We also write $X_{k,n}^x$ ($n \ge k$) for the reflected  
walk starting at time $k$ at $x$, so that $X_n^x = X_{0,n}^x$.
Note that we always have 
\begin{equation}\label{eq:decrease}
|X_{k,n+1}^x-X_{k,n+1}^y| \le |X_{k,n}^x-X_{k,n}^y|\,.
\end{equation}

The following is due to \cite{Fe}, \cite{Kn} and \cite{Le}. 

\begin{lem}\label{lem:irreducible} {\rm (a)} 
The reflected random walk with any starting
point is absorbed after finitely many steps by the interval $\Cc$.\\[4pt] 
{\rm (b)} It is topologically irreducible on $\Cc$, that is, 
for every $x \in \Cc$
and open set $B \subset \Cc$, there is $n$ such that 
$p^{(n)}(x,B) = \Prob[X_n \in B \mid X_0=x] > 0\,.$\\[4pt] 
{\rm (c)} The measure $\nu$ on $\Cc$ given by
$$
\nu(dx) = \bigl(1-F_{\mu}(x)\bigr)\,dx\,,
$$
where $dx$ is Lebesgue measure, is an invariant measure for $p(\cdot,\cdot)$. 
\end{lem}

From \eqref{eq:decrease}, one deduces the following.

\begin{lem}\label{lem:transient}
$\qquad\qquad \quad\qquad\displaystyle \Prob[X_n^x \to \infty] \in \{0\,,1\}\,,$\\
and the value is the same for each starting point $x$.
\end{lem}

\begin{proof} By \eqref{eq:decrease}, the event $[X_n^x \to \infty]$
is in the tail $\sigma$-algebra of  the $(Y_n)$.
\end{proof}

If $\Prob[X_n^x \to \infty]=1$, then we call the reflected random walk
\emph{transient}. 

We now state two important results that were proved 
in \cite{Le} in the case when $\Ex(Y_1) < \infty$, and in the
general case in \cite{Be}. 

\begin{pro}\label{pro:contractive} 
In the non-lattice case, reflected random walk is \emph{locally contractive:}
for every bounded interval $I \subset \Cc$ and all $x, y \in \Cc$,
$$
\lim_{n\to \infty} \uno_I(X_n^x)\,|X_n^x-X_n^y| = 0 \quad \text{almost surely.} 
$$
If $\Prob[X_n^x \to \infty]=0$, then one even has 
$$
\lim_{n\to \infty} |X_n^x-X_n^y| = 0 \quad \text{almost surely.} 
$$
\end{pro}

Of course, also reflected random walk started at time $k$ is locally contractive
for each $k \ge 0$. 
The proof of Proposition \ref{pro:contractive} is outlined in the Appendix.

For $\omega \in \Omega$, let $\Ll^x(\omega)$ be the set of (finite) accumulation 
points of the sequence $X_n^x(\omega)$. In the transient case, $\Ll^x(\cdot)$
is almost surely empty. Otherwise, contractivity implies that there is a set
$\Ll \subset \Cc$, the \emph{attractor} of the process, 
such that 
\begin{equation}\label{eq:attractor}
\Prob[ \Ll^x(\cdot) = \Ll \;\ \text{for all}\;\ x \in \Cc] = 1\,.
\end{equation}
Thus, for any $x\in \Cc$, every open set that interesects $\Ll$ is visited 
infinitely often by $(X_n^x)$ with probability~$1$. 
In other words, the attractor $\Ll$
is \emph{topologically recurrent}, so that it is justified to call the
random walk \emph{recurrent} when $\Prob[X_n^x \to \infty]=0$. 

Proposition \ref{pro:contractive} has the following important consequence,
see the Appendix.

\begin{pro}\label{pro:uniqueness}
In the recurrent case, $\Ll = \Cc\,$, and 
the measure $\nu$ defined in Lemma \ref{lem:irreducible}.c is the unique 
invariant Radon measure for $p(\cdot,\cdot)$ up to multiplication with 
constants. 
\end{pro}

Thus, we have topological recurrence on the whole of $\Cc\,$.
Now, $\nu$ is invariant even in the transient case. If $\Ex(Y_1) < \infty$ then
$\nu(\Cc) < \infty\,$, and we have recurrence by \cite{Le}. As in the
lattice case, we want to extend this recurrence criterion.
Here is the continuous analogue of Theorem \ref{thm:ref-inv} regarding
the process of reflections of \S \ref{sect:reflection}, with 
a rather similar proof.

\begin{thm}\label{thm:ref-inv-cont} In the non-lattice case,
the measure $\rho$ on $\Cc$, given by
$$
\rho(dx)= 
\left( \int_{[0\,,\,\infty)} \mu\bigl((x\,,\,x+y]\bigr) \,\mu(dy)\right) dx
$$
is an invariant measure for the process of reflections $(R_k)$. It is unique
(up to multiplication with constants), if the measure $\nu$ is the unique
invariant measure for the reflected random walk
(up to multiplication with constants).
\end{thm} 

\begin{proof}
We use once more the convolution identity $\Aa * \Uu = \Uu * \Aa = \de_0\,$, 
where $\Aa = \de_0 - \mu$. For any Radon measure $\Mm$ on $\R$, we 
denote by $\check\Mm$ its reflection: 
$\check\Mm(B) = \Mm(-B)$ for Borel sets $B \subset \R$.
We write again $H(x) = 1-F_{\mu}(x)$ for the density of $\nu$ with respect to 
Lebesgue measure, and $h(x)$ for the density of $\rho$. Then
$$
h(x) = \int_{[0\,,\,\infty)} \bigl( H(x) - H(x+y) \bigr)\,\mu(dy)
= H(x) - \check\mu*H(x) = \check\Aa*H(x)\,,
$$
that is, $\rho = \check\Aa*\nu$.  Therefore with the same computations
as above, 
$$
\nu(B)=\check\Uu*\rho(B) = 
\Ex_{\rho}\left(\sum_{j=0}^{\rb(1)-1}\uno_B(X_j)\right)\,,
$$ 
where of course we intend $\Ex_{\rho} = \int \Ex_w(\cdot)\,\rho(dw)$.
Now invariance of $\nu$ for the reflected random walk implies
invariance of $\rho$ for the process of reflections precisely as
in the proof of Theorem \ref{thm:ref-inv}.

For proving uniqueness, let $\bar\rho $ be an invariant (Radon) measure
for $(R_k)$. Once we can prove that the convolution
$\bar\nu  = \check\Uu * \bar\rho $, restricted to $[0\,,\,\infty)$, defines a 
Radon measure (i.e., is finite on compact sets), we can proceed as before: 
$\bar\nu $ is invariant for $(X_n)$, whence $\bar\nu  = c \cdot \nu$ for some $c > 0$, and
$\bar\rho  = \check\Aa * \bar\nu  = c \cdot \check\Aa * \nu = c \cdot \rho$.

If $N < \infty$ then $\bar\rho $ has finite mass, since it
must be concentrated on $[0\,,\,N]$ by Lemma \ref{lem:irreducible}.
Let $\Uu_N$ be the restriction of $\Uu$ to $[0\,,\,N]$. It is also
a finite measure, and on $[0\,,\,N]$, we have 
$\check\Uu * \bar\rho  =\check\Uu_N * \bar\rho $, which is again finite.

Now suppose that $N=\infty$. Let $a > 0$. Then $\supp(\mu)$ contains an
element $M > a$. Choose $b$ such that $a<b<M$. 
Now let $f$ be a compactly supported, continuous function on $\R$, supported
within $[0\,,\,\infty)$, such that $f \equiv 1$ on $[M-b\,,\,M+b]$.  
Then the convolution 
$$
\mu * \check f(x) = \int_{[0\,,\,\infty)} f(v-x)\,\mu(dv) 
$$
defines a continuous function. 

If $x \in [0\,,\,a]$ then $f(v-x)=1$
for all $v \in [M-b+x\,,\,M+b+x] \supset [M-b+a\,,\,M+b]$.
Therefore 
$$
\mu * \check f(x) \ge \mu\bigl([M-b+a\,,\,M+b]\bigr) > 0 \quad
\text{for each} \quad x \in [0\,,\,a]\,.
$$
Using this, the invariance of $\bar\rho $ for $(R_k)$,
the formula of Lemma \ref{lem:reflection}, and Fubini's theorem, we
now compute the finite number
$$
\begin{aligned} 
\int_{[0\,,\,\infty)} f(x)\,\bar\rho (dx) 
&= \int_{[0\,,\,\infty)} \int_{[0\,,\,\infty)} f(y)\;q(x,dy)\;\;\bar\rho (dx) \\
&= \int_{[0\,,\,\infty)} \int_{[0\,,\,\infty)} \int_{[0\,,\,x)}
f(y+w-x)\;\;\Uu(dw)\;\mu(dy)\;\bar\rho (dx) \\
&= \int_{[0\,,\,\infty)} \int_{[0\,,\,\infty)} \int_{[w\,,\,\infty)} 
f(y+w-x) \;\;\bar\rho (dx)\;\Uu(dw)\;\mu(dy) \\
&= \int_{[0\,,\,\infty)} \int_{[0\,,\,\infty)} \int_{[0\,,\,\infty)} 
f(y-x) \;\;[\de_{-w}*\bar\rho](dx)\;\Uu(dw)\;\mu(dy) \\
&= \int_{[0\,,\,\infty)} \int_{[0\,,\,\infty)} 
\mu * \check f(x) \;\;[\de_{-w}*\bar\rho](dx)\;\Uu(dw)\\
&= \int_{[0\,,\,\infty)} \mu * \check f(x) \;\; [\check\Uu*\bar\rho](dx)\\
&\ge \mu\bigl([M-b+a\,,\,M+b]\bigr)\;
\check\Uu*\bar\rho  \bigl( [0\,,a] \bigr)\,. 
\end{aligned}
$$
Therefore $\check\Uu*\bar\rho  \bigl( [0\,,a] \bigr)$ is finite for each
$a>0$.
\end{proof} 

The following is now obtained precisely as in the lattice case.

\begin{cor}\label{cor:finitemass}
The invariant measure $\rho$ of the process of reflections has finite mass
if and only if
\begin{equation}\label{eq:nonl-quadratic}
\int_{[0\,,\,\infty)} \bigl(1 - F_{\mu}(x)\bigr)^2 \,dx < \infty\,.
\end{equation}
This holds, in particular, when 
$
\Ex\bigl(\sqrt{Y_1}\bigr) = \int_{[0\,,\,\infty)} \sqrt{x}\,\mu(dx) < \infty\,.
$
\end{cor}

We now want to deduce recurrence of reflected random walk. This is not as
straightforward as in the case of Markov chains with a denumerable state
space. 

\begin{pro}\label{pro:tail} 
Let $J=(a,b) \subset \Cc$ be a bounded, open interval. Then, setting 
$J(\ep)=(a+\ep\,,\,b-\ep)$,
$$
\Prob\left[\exists\ \ep > 0\,:\; \sum_{n=0}^{\infty} \uno_{J(\ep)}(X_n^x) = \infty \right]
\in \{0, 1 \}\,.
$$
\end{pro}

\begin{proof} 
Each of the countably many events
$$
\left[\lim_{n\to \infty} \uno_{[0\,,\,m]}(X_{k,n}^{\bar x})\,
  |X_{k,n}^{\bar x}-X_{k,n}^{\bar y}| = 0\right] \subset \Omega\,,
$$
where $\bar x,\bar y \in \Cc$ are rational and $k, m \in \N_0$,
has probability $1$. Let $\Omega_0$ be their intersection, so that
$\Prob(\Omega_0)=1$. 
Consider the event
$$
A_J^x = \Omega_0 \cap \bigcup_{0 < \ep < (b-a)/2} B_{J(\ep)}^x\,,\quad\text{where}
\quad B_{J(\ep)}^x = \left[ \sum_{n=0}^{\infty} \uno_{J(\ep)}(X_n^x) = \infty \right]\,. 
$$
We claim that $A_J^x$ does not depend on $x$. Let $y \in \Cc$. 
If $\omega \in A_J^x$ then there is $\ep \in (0\,,\,\frac{b-a}{2})$ 
such that $\omega \in B_{J(\ep)}^x$. 
There are rational numbers $\bar x, \bar y \in \Cc$ such that $|x - \bar x| < \ep/4$
and $|y - \bar y| < \ep/4$. Since $\omega \in \Omega_0$, we have
$$
\uno_{J}\bigl(X_{n}^{\bar x}(\omega)\bigr)\,
\bigl|X_{n}^{\bar x}(\omega)-X_{n}^{\bar y}(\omega)| < \ep/4
$$
for all sufficiently large $n$. Since $|X_n^x - X_n^{\bar x}| \le |x - \bar x|$
and $|X_n^y - X_n^{\bar y}| \le |y - \bar y|$, we get that $X_n^y(\omega) \in J(\ep/4)$
whenever $X_n^x(\omega) \in J(\ep)$. Therefore, $A_J^x \subset A_J^y$, and exchanging
the role of $x$ and $y$, we see that $A_J = A_J^x$ is the same for all $x$.

Now, we claim that $A_J$ is in the tail $\sigma$-algebra of the $(Y_n)_{n\ge 1}$.
Let $\omega \in A_J$ and $\bar \omega \in \Omega$ such that for some $k \in \N$,
$Y_n(\bar\omega) = Y_n(\omega)$ for $n > k$. Then clearly 
$\bar\omega \in \Omega_0$. Set $u=Y_k(\omega)$ and $v=Y_k(\bar\omega)$.
Then we have $X_n^x(\omega) = X_{k,n}^u(\omega)$ and $X_n^x(\bar\omega) = 
X_{k,n}^v(\omega)$ for all $n \ge k$. Now the same ``$\ep/4$''-argument as above
implies that $\bar\omega \in A_J$.

Therefore $\Prob( A_J ) \in \{0\,,\,1\}$ by the 0-1 law of Kolmogorov.
\end{proof}

\begin{thm}\label{thm:nonl-recurrent}
Suppose that the ``quadratic tail'' condition \eqref{eq:nonl-quadratic}
holds. Then, for every starting point $x > 0$, the reflected random walk 
$(X_n^x)$ is topologically recurrent: for every bounded, open interval 
$J \subset \Cc$,
$$
\Prob\left[ \sum_{n=0}^{\infty} \uno_J(X_n^x) = \infty \right] =1\,.
$$
If in addition $\Ex(Y_1) = \int_{[0\,,\,\infty)} x\,\mu(dx) < \infty$,
then $(X_n^x)$ is positive recurrent, while it is null recurrent when 
$\Ex(Y_1)=\infty$.
\end{thm}

\begin{proof}
We write $(R_n^x)$ for the process of reflections starting at $x \in \Cc$, and
define
$$
M_n^x = \frac{1}{n} \sum_{k=0}^{n-1} \uno_{J(\ep)}(R_k^x)\AND 
M^x = \limsup_{n \to\infty} M_n^x\,,
$$
where $\ep > 0$ is chosen such that $J(\ep)$ is non-empty.
The measure $\rho$ of Theorem \ref{thm:ref-inv-cont}
is supported by the whole of $\Cc$, and $\rho(\Cc) < \infty$ by
assumption. We have
$$
\int_{[0\,,\,\infty)} \int_{\Omega} M_n^x \;d\,\Prob\;\rho(dx) 
= \rho\bigl(J(\ep)\bigr)
$$
Since $\rho(\Cc) < \infty$ by assumption and $0 \le M_n \le 1$, we may apply the
``$\limsup$''-variant of the Lemma of Fatou to obtain
$$
\int_{[0\,,\,\infty)} \int_{\Omega} M^x \;d\,\Prob\;\rho(dx) 
\ge \rho\bigl(J(\ep)\bigr)\,.
$$
Therefore there must be $x \in \Cc$ such that
$$
\int_{\Omega} M^x \;d\,\Prob 
\ge \frac{3\rho\bigl(J(\ep)\bigr)}{4\rho(\Cc)}\,.
$$
Consequently, 
$$
0 < \Prob[M_x \ge c] 
\le \Prob\left[\sum_{n=0}^{\infty} \uno_{J(\ep)}(X_n^x) = \infty\right]\,, 
\quad\text{where}\quad 
c = \frac{\rho\bigl(J(\ep)\bigr)}{2\rho(\Cc)} > 0\,.
$$
Proposition \ref{pro:tail} now yields that
$$
\Prob\left[\exists\ \ep > 0\,:\; \sum_{n=0}^{\infty} \uno_{J(\ep)}(X_n^x) 
= \infty \right] = 1\,,
$$
and the result follows.
\end{proof}

Note that we should be careful in stating that the process of reflections
itself is topologically recurrent on $\Cc$ when it has a finite invariant
Radon measure. Indeed, it is by no means clear that it inherits 
local contractivity, or even the property to be Fellerian, 
from reflected random walk.

\section{General reflected random walk}\label{Poisson}

In this section, we drop the restriction that the random variables $Y_n$
are non-negative. Thus, the ``ordinary'' random walk 
$S_n = Y_1 + \cdots + Y_n$ may visit the positive
as well as the negative half-axis. 
Again, $\mu$ will denote the distribution of each of the $Y_n$. 
In the lattice case, we suppose without loss of generality that 
$\supp(\mu) \subset \Z$, and that the group generated by $\supp(\mu)$
is the whole of $\Z$. In the non-lattice case, the closed group
generated by $\supp(\mu)$ is $\R$. 

We start with a simple observation (\cite{Be2} has a more complicated proof).

\begin{lem}\label{lem:symmetric}
If $\mu$ is symmetric, then reflected random walk is (topologically) recurrent 
if and only if the random walk $S_n$ is recurrent.
\end{lem}

\begin{proof} If $\mu$ is symmetric, then also $|S_n|$ is a Markov chain.
Indeed, for a Borel set $B \subset [0\,,\,\infty)$, 
$$
\begin{aligned}
\Prob[\,|S_{n+1}| \in B \mid S_n=x] 
&= \mu(-x+B) + \mu(-x-B) - \mu(-x)\, \de_0(B)\\
&=\Prob[\,|S_{n+1}| \in B \mid S_n=-x]\,,
\end{aligned}
$$ 
and we see that $|S_n|$ has the same transition probabilities as the 
reflected random walk governed by $\mu$. 
\end{proof}   

Recall the classical result that when $\Ex(|Y_1|) <\infty$ and
$\Ex(Y_1)=0$ then $S_n$ is recurrent; see {\sc Chung and Fuchs~\cite{ChFu}}. 
So if $\mu$ is 
symmetric and has finite first moment then reflected random walk is recurrent. 

In general, we should exclude that $S_n \to -\infty$, since in that 
case there are only finitely many
reflections, and reflected random walk tends to $+\infty$ almost surely.

Let $Y_n^+ = \max \{Y_n, 0 \}$ and $Y_n^- = \max \{-Y_n, 0\}$, so that
$Y_n = Y_n^+ - Y_n^-$. The following is well-known.
\begin{lem}\label{lem:infty} 
If \emph{(a)} $\;\Ex(Y_1^-) < \Ex(Y_1^+) \le \infty\,$, or if
\emph{(b)} $\;0 < \Ex(Y_1^-) = \Ex(Y_1^+) < \infty\,$, then
$\limsup S_n = \infty\,$ almost surely, so that there are infinitely many
reflections.
\end{lem}

We note that Proposition \ref{pro:contractive} is also valid here,
since its proof (see the Appendix) does not require non-negativity
of $Y_n$. Also, when the $Y_n$ may assume both positive and negative 
values with positive probability, then the essential class, resp. classes,
on which reflected random walk evolves is/are unbounded. In the non-lattice
case this is $\Cc = [0\,,\,\infty)$, and $X_n^x$ is locally contractive.

In the sequel, we assume that $\limsup S_n = \infty$ almost surely.
Then the (non-strictly) ascending  \emph{ladder epochs}
$$
\lb(0)  = 0\,,\quad  \lb(k+1) = \inf \{ n > \lb(k) : S_n \ge S_{\lb(k)} \}
$$ 
are all almost surely finite, and the random variables $\lb(k+1) - \lb(k)$ 
are i.i.d.
We can consider the \emph{embedded random walk} $S_{\lb(k)}\,$, $k \ge 0$, 
which tends to $\infty$ almost surely. Its increments 
$\overline Y_k = S_{\lb(k)} - S_{\lb(k-1)}\,$, $k \ge 1$, are i.i.d. non-negative
random variables with distribution denoted $\overline{\mu}$. 
Furthermore, if $\overline{\!X}_k^x$ denotes the reflected random
walk associated with the sequence $(\overline Y_k)$, while $X_n^x$ is our original
reflected random walk associated with $(Y_n)$, then 
$$
\overline{\!X}_k^x = X_{\lb(k)}^x\,,
$$ 
since no reflection can occur between times $\lb(k)$ and $\lb(k+1)$.

\begin{lem}[\cite{Be}]\label{lem:recequ}
The embedded reflected random walk $\overline{\!X}_k^x$ is recurrent if and 
only the original reflected random walk is recurrent.
\end{lem} 

\begin{proof} Since both processes are locally contractive, each of the
two processes is transient if and only if it tends to $+\infty$ almost surely:
in the lattice case this is clear, and in the non-lattice case it follows 
from local contractivity. 
If $\lim_n X_n^x = \infty$ then clearly also $\lim_k X_{\lb(k)}^x = \infty$ a.s.
Conversely, suppose that $\lim_k \overline{\!X}_k^x \to \infty$ a.s.
If $\lb(k) \le n < \lb(k+1)$ then $X_n^x \ge X_{\lb(k)}^x$.
(Here, $k$ is random, depending on $n$ and $\omega \in \Omega$, and when
$n \to \infty$ then $k \to \infty$ a.s.) Therefore, also 
$\lim_n X_n^x = \infty$ a.s.
\end{proof}

As long as $\limsup S_n = \infty\,$, we can consider the reflection 
times as in \eqref{eq:reflectiontime} for the case of non-negative $Y_n$.
The observation that there is no reflection between 
times $\lb(k)$ and $\lb(k+1)$ yields the following.

\begin{lem}\label{lem:sametimes}
The reflection times for $(X_n^x)$ and $(\overline{\!X}_k^x)$ are the same,
so that reflected random walk and embedded reflected random walk have the
same process of reflections. In particular, if the latter has a finite invariant
measure, resp., if it is non-transient, then $(X_n^x)$ is (topologically)
recurrent on its essential class(es).
\end{lem}

We can now deduce the following.

\begin{thm}\label{thm:sqrt}
Reflected random walk  $(X_n^x)$ is (topologically) recurrent on its essential 
class(es), if
\begin{quote} 
\emph{(a)} $\;\Ex(Y_1^-) < \Ex(Y_1^+)$ and 
$\Ex\bigl(\sqrt{Y_1^+}\,\bigr) < \infty\,,$ or if\\
\emph{(b)} $\;0 < \Ex(Y_1^-) = \Ex(Y_1^+)$ and 
$\Ex\Bigl(\sqrt{Y_1^+}^{\,3}\Bigr) < \infty\,$. 
\end{quote}
\end{thm}

\begin{proof}
We show that in each case the assumptions imply that 
$\Ex\bigl(\sqrt{\,\overline{Y}_1}\bigr) < \infty$.
Then we can apply Lemma \ref{lem:momenthalf}, resp. Corollary 
\ref{cor:finitemass} to deduce recurrence of $(\overline{\!X}_k^x)$.
This in turn yields recurrence of $(X_n^x)$ by Lemma \ref{lem:sametimes}.

\smallskip

(a) Under the first set of assumptions,
$$
\begin{aligned}
\Ex\Bigl(\sqrt{\overline{Y}_1}\Bigr) 
&= \Ex\Bigl(\sqrt{Y_1+\ldots+Y_{\lb(1)}^{\,} }\,\Bigr)
\le \Ex\Bigl(\sqrt{Y_1^+ +\ldots+Y_{\lb(1)}^+}\,\Bigr)\\
&\le \Ex\Bigl(\sqrt{Y_1^+}+\ldots+\sqrt{Y_{\lb(1)}^+}\,\Bigr) 
= \Ex\Bigl(\sqrt{Y_1^+}\,\Bigr) \cdot\Ex\bigl(\lb(1)\bigr) 
\end{aligned}
$$
by Wald's identity. Thus, we now are left with proving
$\Ex\bigl(\lb(1)\bigr) < \infty\,$. If $\Ex(Y_1^+) < \infty$, then
$\Ex(|Y_1|) < \infty$ and $\Ex(Y_1) > 0$ by assumption,
and in this case it is well known that $\Ex\bigl(\lb(1)\bigr) < \infty\,$;
see e.g. \cite[Thm. 2 in \S XII.2, p. 396-397]{Fe}. 
If  $\Ex(Y_1^+) = \infty$ then there is $M > 0$ such that
$Y_n^{(M)} = \min\{ Y_n\,, M\}$ (which has finite first moment)
satisfies $\Ex(Y_n^{(M)}) = \Ex(Y_1^{(M)})> 0\,$. The first increasing
ladder epoch $\lb^{(M)}(1)$ associated with 
$S_n^{(M)} = Y_1^{(M)} + \ldots + Y_n^{(M)}$ has finite expectation by
what we just said, and $\lb(1) \le\lb^{(M)}(1)$. Thus, $\lb(1)$
is integrable.

\smallskip

(b) If the $Y_n$ are centered, non-zero and 
$\Ex\bigl((Y_1^+)^{1+a}\bigr)< \infty\,,$ where $a > 0$, then 
$\Ex\bigl((\overline{Y}_1)^a \bigr) < \infty\,$, as was shown by 
{\sc Chow and Lai~\cite{ChLa}}. In our case, $a=1/2$.
\end{proof}

In conclusion, we discuss sharpness of the sufficient recurrence conditions
$\Ex\Bigl(\sqrt{Y_1}^{\,3}\Bigr) < \infty$ in the centered case,
resp $\Ex\bigl(\sqrt{Y_1}\bigr) < \infty$ in the case when $Y_1 \ge 0$.

\begin{exa}\label{ex:symm}
Define a symmetric probability measure $\mu$ on $\Z$ by
$$
\mu(0)=0\,,\quad \mu(k)=\mu(-k) = c/k^{1+a}\quad (k \ne 0)\,,
$$
where $a >0$ and $c$ is the proper normalizing constant. Then it is known
that the associated symmetric random walk $S_n$ on $\Z$ is recurrent if and 
only if $a \ge 1$, see {\sc Spitzer~\cite[p. 87]{Sp}}. 
By Lemma \ref{lem:symmetric}, the associated 
reflected random walk is also recurrent, but when $1 \le a \le 3/2$ then
condition (b) of Theorem \ref{thm:sqrt} does not hold. 
\end{exa} 

Nevertheless, we can also show that in general, the sufficient condition 
$\Ex\Bigl(\sqrt{\,\overline{Y}_1\,}\Bigr) < \infty$ for
recurrence of reflected random walk with non-negative increments 
$\overline{Y}_n$ is very close to being sharp. (We write $\overline{Y}_n$
because we shall represent this as an embedded random walk in the next 
example.)

\begin{pro}\label{pro:nsymm}
Let $\mu_0$ be a probability measure on $\N_0$ such that
$\mu_0(n) \ge \mu_0(n+1)$ for all $n \ge 0$ and 
$$
\mu_0(n) \sim c \,(\log n)^b \big/ n^{3/2}\,,
\quad\text{as}\;n\to \infty\,,
$$
where $b > 1/2$ and $c > 0$.
Then the associated reflected random walk on $\N_0$ is transient.
\end{pro}

Note that $\mu_0$  has finite moment of order $\frac12 - \ep$ for 
every $\ep > 0$, while the moment of order $\frac12$ is infinite. 

\smallskip

The proof needs some preparation. Let $(Y_n)$ be i.i.d. random variables
with values in $\Z$ that have finite first moment and are non-constant
and centered, and let $\mu$ be their common distribution. 
The first \emph{strictly ascending} and 
\emph{strictly descending ladder epochs} of the random walk
$S_n=Y_1 + \ldots + Y_n$ are 
$$
\tb_+(1)  = \inf \{ n > 0 : S_n > 0 \} \AND
\tb_-(1)  = \inf \{ n > 0 : S_n < 0 \}\,,
$$
respectively. They are almost surely finite. Let $\mu_+$ be the distribution
of $S_{\tb_+(1)}$ and $\mu_-$ the distribution of $S_{\tb_-(1)}$, and --
as above -- $\overline{\mu}$ the distribution of $\overline{Y}_1 = S_{\lb(1)}$.
We denote the characteristic function associated with any probability measure
$\sigma$ on $\R$ by $\wh\sigma(t)\,$, $t \in \R$. Then, following 
{\sc Feller~\cite[(3.11) in \S XII.3]{Fe}}, \emph{Wiener-Hopf-factorization}
tells us that
\begin{equation}\label{eq:WH}
\begin{gathered}
\mu = \overline{\mu} + \mu_- - \overline{\mu}*\mu_-
\;\AND\;
\overline{\mu} = u\cdot \de_0 + (1-u)\cdot\mu_+\ ,\\
\text{where}\quad
u = \overline{\mu}(0) = 
\sum_{n=1}^{\infty} \Prob[S_1 < 0\,, \ldots, S_{n-1} < 0\,,\; S_n =0] < 1\,. 
\end{gathered}
\end{equation} 
(Recall that $*$ is convolution.)
 
\begin{lem}\label{lem:ladder}
Let $\mu_0$ be a probability measure on $\N_0$ such that
$\mu_0(n) \ge \mu_0(n+1)$ for all $n \ge 0$.
Then there is a symmetric probability measure $\mu$ on $\Z$ such that
that the associated first (non-strictly) ascending ladder random variable
has distribution $\mu_0$.
\end{lem}

\begin{proof}
We decompose $\;\mu_0 = \mu_0(0)\cdot\de_0 + 
\bigl(1-\mu_0(0)\bigr)\cdot \mu_{\times}\,,$ where $\mu_{\times}$ is supported
by $\N$ (i.e., $\mu_{\times}(0)=0$). If $\mu_0$ is the law of the first
strictly ascending  ladder random variable associated with some symmetric
measure $\mu$, then by \eqref{eq:WH} we must have $\mu_- = \check\mu_{\times}$,
the reflection of $\mu_{\times}$ at $0$, and
\begin{equation}\label{eq:mubarmu}
\mu = \mu_0 + \check \mu_{\times}  - \mu_0*\check\mu_{\times}
= \mu_0(0)\cdot \de_0 + \bigl(1-\mu_0(0)\bigr)\cdot
(\mu_{\times} + \check \mu_{\times}  - \mu_{\times} * \check\mu_{\times})\,.
\end{equation}
We \emph{define} $\mu$ in this way.
The monotonicity assumption on $\mu_0$ implies that $\mu$ is a probability 
measure: indeed, it is straightforward to
show that $\mu(k) \ge 0$ for each $k \in \Z$.

The measure $\mu$ of \eqref{eq:mubarmu} is non-degenerate and symmetric. 
If it induces 
a recurrent random walk $(S_n)$, then the ascending and descending ladder
epochs are a.s. finite. If $(S_n)$ is transient, then $|S_n| \to \infty$ 
almost surely, but it cannot be $\Prob[S_n \to \infty] > 0$ since
in that case this probaility had to be 1 Kolmogorov's 0-1-law, while 
symmetry would yield 
$\Prob[S_n \to -\infty] = \Prob[S_n \to \infty] \le 1/2$. 
Therefore $\liminf S_n = -\infty$ and $\limsup S_n = +\infty$ almost
surely, a well-known fact, see e.g. \cite[Thm. 1 in \S XII.2, p. 395]{Fe}. 
Consequently, the ascending and descending ladder 
epochs are again a.s. finite. Therefore the probability measures 
$\mu_+$ and $\mu_-=\check\mu_+$ (the laws of $S_{\tb_\pm(1)}$) are well defined.
By the uniqueness theorem of Wiener-Hopf-factorization
\cite[Thm. 1 in \S XII.3, p. 401]{Fe}, it follows that $\mu_- = \check\mu_{\times}$
and that the distribution of the first (non-strictly) ascending ladder random 
variable is $\overline{\mu} = \mu_0$.
\end{proof}

\begin{proof}[Proof of Proposition \ref{pro:nsymm}]  
Let $\mu$ be the symmetric measure associated with $\mu_0$
according to \eqref{eq:mubarmu} in Lemma \ref{lem:ladder}. 
Then its characteristic function
$\wh\mu(t)$, given by \eqref{eq:WH}, is non-negative real. A well-known
criterion says that the random walk $S_n$ associated with $\mu$ is transient
if and only if (the real part of) $1\big/\bigl(1-\wh\mu(t)\bigr)$ is integrable
in a neighbourhood of $0$. 
Returning to $\overline{\mu}$, it is a standard exercise 
(see  \cite[Ex. 12 in Ch. XVII, Section 12]{Fe})
to show that there is $A \in \C\,$, $A \ne 0$
such that its characteristic function satisfies 
$$
\wh{\overline{\mu}}(t)= 1+A\,\sqrt{t}\,(\log t)^b\,\bigl(1+o(t)\bigr)
\quad \text{as} \; t \to 0\,.
$$
By \eqref{eq:WH},
$$
1 - \wh\mu(t) = 
(1-u)\bigl(1 - \wh{\mu_+}(t)\bigr)\bigl(1 - \wh{\mu_-}(t)\bigr)\,. 
$$
We deduce    
$$
\wh{\mu}(t)=  1 + \bigl(1-\mu_0(0)\bigr)\,|A|^2 \, t \,(\log t)^{2b}\,
\bigl(1+o(t)\bigr)
\quad \text{as} \; t \to 0\,.
$$
The function $1\big/\bigl(1-\wh\mu(t)\bigr)$ is integrable near $0$. 
By Lemma \ref{lem:symmetric}, 
the associated reflected random walk is transient. But then also 
the embedded reflected random walk associated with $S_{\lb(n)}$ is
transient by Lemma \ref{lem:recequ}. This is the reflected random walk
governed by~$\overline{\mu}$.
\end{proof}

\section{Appendix: local contractivity}\label{sect:appendix}
  
Here, we come back to propositions \ref{pro:contractive} and
\ref{pro:uniqueness}. They arise as special cases of two main results 
in the PhD thesis of {\sc Benda~\cite{Be}} and of the contents of the two
papers \cite{Be1} and \cite{Be2}, which were accepted for
publication but remained unpublished. For this reason, we give an outline,
resp. published references for
their proofs. In \cite{Be1}, this is placed in the following more general
context. Let $(\Xx,d)$ be a proper metric space (i.e., closed balls are
compact), and let $\Gp$ be the monoid of all continuous 
mappings $\Xx \to \Xx$.
It carries the topology of uniform convergence on compact sets.
Now let $\wt\mu$ be a regular probability measure on $\Gp$, and let 
$(F_n)_{n \ge 1}$ be a sequence of i.i.d. $\Gp$-valued random variables
(functions) with common distribution $\wt\mu$. 
The measure $\wt\mu$ gives rise to the stochastic iterated function system
(SFS) $X_n^x$ defined by
\begin{equation}\label{eq:SFS}
X_0^x = x \in \Xx\,,\AND X_n^x = F_n(X_{n-1}^x)\,,\quad n \ge 1\,.
\end{equation}
In the setting of the above Sections
\ref{sect:reflection}--\ref{sect:non-lattice}, we have $\Xx = [0\,,\,\infty)$
with the standard distance, and $F_n(x) = |x-Y_n|$, so that the measure
$\wt\mu$ is the image of the distribution
$\mu$ of the $Y_n$ in \S \ref{sect:reflection} under the mapping
$[0\,,\,\infty) \to \Gp\,$, $y \mapsto g_y$, where $g_y(x) = |x-y|$.

\begin{dfn}\label{def:loccont}
The SFS is called \emph{locally contractive,} if for 
all $x \in \Xx$ and every compact $K \subset \Xx$, 
$$
\uno_K(X_n^x) \cdot \sup_{y \in K} d(X_n^x,X_n^y) \to 0 \quad\text{almost surely, as}\; 
n \to \infty\,.
$$
\end{dfn}

This notion was first introduced by {\sc Babillot, Bougerol and
Elie~\cite{BaBoEl}} and was later exploited systematically by {\sc Benda},
who (in personal comunication) also gives credit to unpublished work
of his late PhD advisor {\sc Kellerer}, compare with the posthumous
publication \cite{Ke}.

Using Kolomogorov's 0-1 law (and properness of $\Xx$), one   gets a 
general variant of Lemma~\ref{lem:transient}.
\begin{lem}\label{lem:transient2} For a locally contractive SFS of 
contractions,
$$
\begin{aligned}
\text{either}\quad &\Prob[d(X_n^x,x) \to \infty]=0 
\quad\text{for all}\;\ x \in \Xx\,,\\ 
\text{or}\qquad\;\; &\Prob[d(X_n^x,x) \to \infty]=1 
\quad\text{for all}\;\ x \in \Xx\,.
\end{aligned}
$$
\end{lem} 

\begin{proof} Let $B(r)\,$, $r \in \N$ be the open balls
in $\Xx$ with radius $r$ and fixed center $o \in \Xx$. It has compact
closure by properness of $\Xx$. Consider 
\begin{equation}\label{eq:Xmn}
X_{m,n}^x = F_n \circ F_{n-1} \circ \ldots \circ F_{m+1}(x)
\end{equation}
for $n > m$, so that $X_n^x =X_{0,n}^x\,$. Then local contractivity
implies that for each $x \in \Xx$, we have $\Prob(\Omega_0)=1$ for the
event $\Omega_0$ consisting of all $\omega \in \Omega$ with
\begin{equation}\label{eq:Om0}
\lim_{n\to \infty} \uno_{B(r)}\bigl(X_n^x(\omega)\bigr) \cdot 
\sup_{y \in B(r)} d\bigl(X_{m,n}^x(\omega),X_{m,n}^y(\omega)\bigr) = 0
\quad\text{for each}\; r \in \N\,,\; m \in \N_0\,. 
\end{equation}
Clearly, $\Omega_0$ is invariant with respect to the shift of the sequence 
$(F_n)$. 
 
Now let  $\omega \in \Omega_0$ be such that the sequence 
$\bigl( X_n^x(\omega)\bigr)_{n \ge 0}$ accumulates at some $w \in \Xx$.
Fix $m$ and set $v = X_m^x(\omega)$. Then also 
$\bigl(X_{m,n}^{v}(\omega)\bigr)_{n \ge m}$
accumulates at $w$. Now let $y \in \Xx$ be arbitrary. Then there is $r$ such
that $v, w, y \in B(r)$. Therefore also 
$\bigl( X_{m,n}^y(\omega)\bigr)_{n \ge m}$ accumulates at $w$. 
In particular, the fact that $\bigl( X_n^x(\omega)\bigr)_{n \ge 0}$ accumulates
at some point does not depend on the initial trajectory, i.e., on the
specific realization of $F_1, \dots, F_m$. We infer that the set
$$
\bigl\{ \omega \in \Omega_0 : \bigl( X_n^x(\omega)\bigr)_{n \ge 0}
\; \text{accumulates in} \; \Xx \bigr\}
$$ 
is a tail event of $(F_n)_{n \ge 1}$. On its complement in $\Omega_0$,
we have $d(X_n^x,x) \to \infty\,$.
\end{proof}

If $d(X_n^x,x) \to \infty$ almost surely, then we call the SFS \emph{transient.}
What has been said about the attractor in \eqref{eq:attractor} for\; 
reflected random walk is true in general. For $\omega \in \Omega$, let
$\Ll^x(\omega)$ be the set of accumulation points of $\bigl( X_n^x(\omega)\bigr)$
in $\Xx$. A straightforward extension of the argument used in the last
proof (using again properness of $\Xx$) yields the following.

\begin{lem}\label{lem:attract} 
For any non-transient, locally contractive SFS, there is a set 
$\Ll \subset \Xx$ -- the \emph{attractor} -- such that
$$
\Prob[ \Ll^x(\cdot) = \Ll \;\ \text{for all}\;\ x \in \Cc] = 1\,,
$$
\end{lem}
Thus, $(X_n^x)$ is \emph{(topologically) recurrent} on $\Ll$
when $\Prob[d(X_n^x,x) \to \infty]=0$. 
 
\begin{pro}\label{pro:uniquemeasure}
For a recurrent locally contractive SFS, there is a unique invariant Radon
measure $\nu$ on $\Xx$ up to multiplication with constants, and 
$\supp(\nu) = \Ll$.
\end{pro}  
 
This is contained in \cite{Be} and \cite{Be1}. The proof of the existence of
such a measure supported by $\Ll$ is rather straightforward, compare with
the old survey by Foguel~\cite{Fo}. (One first constructs an excessive
measure supported by $\Ll$ via a ratio limit argument, an then uses recurrence
to obtain that it has to be invariant.) For a proof of uniqueness that
is available in print, see {\sc Brofferio~\cite[Thm. 3]{Br}}, who considers
only SFS of affine mappings, but the argument carries over to general 
locally contractive SFS without changes.
  
Let us now consider a more specific class of SFS:  
within $\Gp$, we consider the closed submonoid $\Lp$ of all 
\emph{contractions} of $\Xx$, i.e., mappings $f: \Xx \to \Xx$ with 
Lipschitz constant $L(f) \le 1$. We suppose that the probability measure
$\wt\mu$ that governs the SFS is supported by $\Lp$, that is, each 
random function $F_n$ of \eqref{eq:SFS} satisfies $L(F_n) \le 1$.
In this case, one does not need local contractivity in order to obtain
Lemma \ref{lem:transient2}; 
this follows directly from properness of $\Xx$ and the inequality
$$
d(X_n^x,X_n^y) \le d(x,y)\,.
$$
Let $\Sf(\wt\mu)$ be the closed sub-semigroup of $\Lp$ generated by 
$\supp(\wt\mu)$. 
The following key result of \cite{Be} is inspired by \cite[Thm. 2.2]{Kn}, where 
reflected random walk with $\Ex(Y_n) < \infty$ is studied. 

\begin{pro}\label{pro:contractive2} If \emph{(i)} the SFS of
contractions is non-transient, 
and \emph{(ii)} the semigroup $\Sf(\wt\mu) \subset \Lp$ contains  a constant 
function, then
$$
D_n(x,y) = d(X_n^x,X_n^y) \to 0 \quad\text{almost surely, as}\; n \to \infty\,.
$$
\end{pro}

\begin{proof}
Since $D_{n+1}(x,y) \le D_n(x,y)$, the limit $D_{\infty}(x,y) = \lim_n D_n(x,y)$
exists and is between $0$ and $d(x,y)$. We set 
$w(x,y) = \Ex\bigl( D_{\infty}(x,y)\bigr)$.
First of all, we claim that
\begin{equation}\label{eq:martingale}
\lim_{m \to \infty} w(X_m^x\,,X_m^y) = D_{\infty}(x,y)\quad\text{almost surely.}
\end{equation}
To see this, consider $X_{m,n}^x$ as in \eqref{eq:Xmn}.  
Then $D_{m,\infty}(x,y) = \lim_n d(X_{m,n}^x,X_{m,n}^y)$ has the same
distribution as $D_\infty(x,y)$, whence 
$\Ex\bigl(D_{m,\infty}(x,y)\bigr) = w(x,y)$. Therefore, we also have
$$
\Ex\bigl(D_{m,\infty}(X_m^x\,,X_m^y) \mid F_1, \ldots, F_m\bigr) = 
w(X_m^x\,,X_m^y)\,.
$$
On the other hand, $D_{m,\infty}(X_m^x\,,X_m^y) = D_{\infty}(x,y)$,
and the bounded martingale
$$
\Bigl(\Ex\bigl(D_{\infty}(x,y) \vert F_1, \ldots, F_m \bigr)
       \Bigr)_{m \ge 1}
$$       
converges  almost surely 
to $D_\infty(x, y)$. The proposed statement \eqref{eq:martingale} follows.

\smallskip

Now let $\ep > 0$ be arbitrary, and fix $x, y \in X$.
We have to show that the event $A = [D_{\infty}(x,y) \ge \ep]$ has probability 
$0$. 

(i) By non-transience, 
$$
\Pr\left( \bigcup_{r \in \N\,} \bigcap_{\,m \in \N\,} \bigcup_{\,n \ge m\,} 
[X_n^x\,,\; X_n^y \in \Bb(r)] \right) = 1\,.
$$
On $A$, we have $D_n(x,y) \ge \ep$ for all $n$. Therefore we need to show
that $\Pr(A_r) = 0$ for each $r \in \N$, where 
$$
A_r = \bigcap_{m \in \N} \bigcup_{n \ge m} 
[X_n^x\,,\; X_n^y \in \Bb(r)\,,\;D_n(x,y) \ge \ep]\,.
$$

(ii) By the second hypothesis, there is $x_0 \in X$ which can be approximated
uniformly on compact sets by functions of the form $f_k \circ ... \circ f_1$,
where $f_j \in \supp(\wt\mu)$. Therefore, given $r$ there is $k \in \N$ 
such that
$$
\Pr(C_{k,r}) > 0\,, \quad \text{where} \quad 
C_{k,r} = \left[\sup_{u \in \Bb(r)} d(X_k^u\,,x_0) \le \ep/4 \right]\,.
$$
On $C_{k,r}$ we have $D_{\infty}(u,v) \le D_k(u,v) \le \ep/2$ for all 
$u,v \in \Bb(r)$. Therefore, setting $\de=\Pr(C_{k,r})\cdot (\ep/2)$,
we have for all $u,v \in \Bb(r)$ with $d(u,v) \ge \ep$ that
$$
\begin{aligned}
w(u,v) &= \Ex\bigl( \uno_{C_{k,r}}\, D_{\infty}(u,v) \bigr) +
\Ex\bigl( \uno_{\Xx \setminus C_{k,r}}\, D_{\infty}(u,v) \bigr) \\
&\le \Pr(C_{k,r})\cdot (\ep/2) + \bigl(1-  \Pr(C_{k,r})\bigr)\cdot d(u,v)
\le d(u,v) - \de\,.
\end{aligned}
$$
We conclude that on $A_r$, there is a (random) sequence $(n_{\ell})$ such that 
$$
w(X_{n_{\ell}}^x\,,X_{n_{\ell}}^y) \le D_{n_{\ell}}(x,y) - \de\,.
$$
Passing to the limit on both sides, we see that \eqref{eq:martingale} is
violated on $A_r$, since $\de > 0$. Therefore $\Pr(A_r) = 0$ for each $r$.
\end{proof}

\begin{cor}\label{cor:loccontractive} If the semigroup $\Sf(\wt\mu) \subset \Lp$
contains  a constant function, then the SFS is locally contractive.
\end{cor}

\begin{proof} In the transient case, $X_n^x$ can visit any compact
$K$ only finitely often, whence $\uno_K(X_n^x)\to 0$ a.s.
In the non-transient case, we use the fact that by properness, $\Xx$
has a dense, countable subset $Y$. Proposition \ref{pro:contractive2}
implies that with probability $1$, we have $\lim_n D_n(x,w)=0$ for all
$w \in Y$. If $K \subset \Xx$ is compact and $\ep > 0$ then there is
a finite $W \subset Y$ such that $d(y, W) < \ep$ for every $y \in K$.
Therefore  
$$
\sup_{y \in K} D_n(x,y) \le 
\underbrace{\max_{w \in W} D_n(x,w)}_{\textstyle{\to 0 \;\text{a.s.}}} 
+ \ep\,,
$$ 
since $D_n(x,y) \le D_n(x,w)+D_n(w,y) \le  D_n(x,w) + d(w,y)$.
\end{proof}

\begin{proof}[Proof of Proposition \ref{pro:contractive}]
Reflected random walk is an SFS of contractions, since $L(g_y)=1$
for the function $g_y(x) = |x-y|$. 
\cite[Prop. 2]{Le}
shows that the constant function $x \mapsto 0$ is contained
in the semigroup $\Sf(\wt\mu)$, where $\mu$ is the law of the
increments $Y_n$ and $\wt\mu$ its image in the semigroup $\Lp$
of contractions of $\Xx = [0\,,\,\infty)$ under the mapping
$y \mapsto g_y$, $g_y(x) = |x-y|$. Note that this statement
and its proof in \cite{Le} are completely deterministic, regarding 
topological properties of the set $\supp(\mu) \subset [0\,,\,\infty)\,$, 
and do not rely on any moment condition.
\end{proof} 

\begin{proof}[Proof of Proposition \ref{pro:uniqueness}]
If reflected random walk is recurrent, then we know from Proposition
\ref{pro:uniquemeasure} that there is a unique invariant Radon measure up to
multiplication with constants, and its support is the attractor $\Ll$.
On the other hand, we already have the invariant measure $\nu$ given in
Lemma \ref{lem:irreducible}.c, and its support is $\Cc$. 
\end{proof}

\end{document}